\documentclass{article}
\usepackage{tikz}
\usepackage{hyperref}
\hypersetup{backref=true,       
	pagebackref=true,               
	hyperindex=true,                
	colorlinks=true,                
	breaklinks=true,                
	urlcolor= black,                
	linkcolor= green,                
	bookmarks=true,                 
	bookmarksopen=false,
	citecolor=black,
	linkcolor=blue,
	filecolor=black,
	citecolor=green,
	linkbordercolor=blue
}

\usepackage[utf8]{inputenc}
\usepackage{lipsum}
\usepackage{graphicx}
\usepackage[toc]{appendix}
\usepackage{pgf,pgfarrows,pgfnodes,pgfautomata,pgfheaps,pgfshade}
\usepackage{rotating}
\usepackage{multirow}
\usepackage{epstopdf}
\usepackage{mathtools}
\usepackage{amsmath,amssymb}
\usepackage{enumitem}
\usepackage{amsthm}
\usepackage{amsfonts}
\usepackage{authblk}
\usepackage{listings}
\usepackage{slashbox}

\usepackage{tikz,xcolor}
\usepackage{tikz-cd}
\usetikzlibrary{shapes.geometric,arrows,positioning,fit,calc,}
\usetikzlibrary{arrows}
\usepackage{imakeidx}
\usepackage{indentfirst}
\usepackage[left=3.5cm, right=3.5cm, top=2.5cm]{geometry}
\usepackage{array}
\usepackage{caption}

\newtheorem{theorem}{Theorem}[section]
\newtheorem*{theorem*}{Theorem}
\newtheorem{lemma}[theorem]{Lemma}
\newtheorem{cor}[theorem]{Corollary}
\newtheorem{prop}[theorem]{Proposition}
\theoremstyle{definition}
\newtheorem{definition}[theorem]{Definition}
\newtheorem{example}[theorem]{Example}
\newtheorem{notation}[theorem]{Notation}
\newtheorem{remark}[theorem]{Remark}

\newcommand\BibTeX{{\rmfamily B\kern-.05em \textsc{i\kern-.025em b}\kern-.08em
		T\kern-.1667em\lower.7ex\hbox{E}\kern-.125emX}}

\title{\textbf{\textsc{CLASSIFICATION OF GROUPLIKE CATEGORIES\thanks{Supported by the Agence Nationale de la Recherche program 3ia Côte d’Azur ANR-19- P3IA-0002, and European Research Council Horizons 2020 grant 670624 (Mai Gehrke’s DuaLL project).\\
				This is a part of my thesis with Carlos Simpson and Khaled Smaili.
}}}}

\author{Najwa Ghannoum \thanks{najwa.ghannoum@univ-cotedazur.fr}}
\affil{\small{Laboratoire J. A. Dieudonn\'e, CNRS, and Universit\'e C\^ote d'Azur, France}}
\affil{Laboratoire LAMA, Universit\'e Libanaise, Lebanon}

\date{}
\begin{document}
	\maketitle
	\begin{abstract}
		In this paper we study grouplike monoids, these are monoids that contain a group to which we add an ordered set of idempotents. We classify finite categories with two objects having grouplike endomorphism monoids, and we give a count of certain categories with grouplike monoids. \\
		
	\end{abstract}
	
	\section{Introduction}
	\textit{Finite categories} are categories with a finite set of objects and a finite set of morphisms. The main purpose of our work in general is the classification of finite categories. We represent finite categories as square matrices such that the entries of the matrix are the number of morphisms between each two objects.
	
	Another way of viewing a finite category, is that we classify monoids of certain orders, then we combine them together in one category (say here for example having two objects), two monoids $\textsf{C}_i, \textsf{C}_j$ are called \textit{connected} if there exists a category with non-empty morphism sets in both directions, where the monoids $\textsf{C}_i$ and $\textsf{C}_j$ are the endomorphism monoids of the two objects. Viewing them as objects, each object in the category is going to correspond to one of the monoid of endomorphisms that we specify. A \textit{groupoid} is a category in which every morphism is invertible.
	
	Starting from this idea, we can represent matrices of categories in a different way. Instead of number of morphisms between objects, we represent them in terms of the endomorphism monoids of the category. For example a category $\textsf{C}$ with two objects $X,Y$ will be represented as follows
	$$
	\left( \begin{array}{cc}
		A & L \\
		R & B
	\end{array} \right)
	$$
	where
	\begin{center}
		$\textsf{C}(X,X) = A$, $\textsf{C}(X,Y) = L$, $\textsf{C}(Y,X) = R$ and $\textsf{C}(Y,Y) = B$
	\end{center}
	such that $A$ and $B$ are monoids with a specific algebraic structure.
	
	From this perspective we can see parts of the category as bimodules, taking only the multiplication tables of $A \cdot L \cdot B$ and $B \cdot R \cdot A$. In some cases, the number of bimodules obtained specifies the number of categories that could be obtained by two bimodules.
	
	We study in this paper a specific type of finite categories called \textit{grouplike} categories. To define grouplike categories, we need to define grouplike monoids first. A \textit{grouplike monoid} is a monoid of the form $G \cup I$ where $G$ is a group and $I = \{e_1,\hdots,e_k\}$ such that
	\begin{equation} \label{idsets} \tag{Ord}
	(e_i \cdot e_j = e_i ~\textrm{and}~ e_j \cdot e_i = e_i) ~\textrm{if and only if} ~ i \leq j.
	\end{equation}
We denote it by $G^{*k}$. A \textit{grouplike category}	is a category where its endomorphism monoids are grouplike.
	
	Counting associative structures has been a problem for years. We introduce some of the previous work in this area:
	
	\textsc{In 2009}, A. Distler and T. Kelsey counted the monoids of orders eight, nine and ten. They weren't able to achieve more than that as the number of semigroups of order 10 was unknown \cite{DK2009}.
	
	\textsc{In 2012}, the number of semigroups of order 10 was known by A. Distler, C. Jefferson, T. Kelsey, and L. Kotthoff \cite{semigroupsorder10}.
	
	\textsc{In 2014},  G. Cruttwell and R. Leblanc introduced the question: \textit{How many categories are there with \textbf{n} morphisms?} It means with a total number of morphisms distributed between objects randomly. They compared the numbers obtained with the number of monoids of order n, which leads to almost the same numbers up to order 10.
	
	\textsc{In 2017}, S. Allouch and C. Simpson counted the categories whose number of morphisms between each two objects is 2. They were able to get an exact count up to order 3, and bounds for a general size \cite{samer2}.
	
	This work is inspired by our results with W. Fussner, T. Jakl and C. Simpson in \cite{AITP2020}. Using the program \textsc{Mace4} \cite{mace4}, we gave a count to the categories associated to the matrix
	$$
	\left( \begin{matrix}
		3 & 3 \\
		3 & 3
	\end{matrix} \right).
	$$
	We classified monoids of order 3 and presented the number of categories between each two monoids. The data obtained from the count showed that the nature of the endomorphism monoids affects the classification in a very important way, which lead us to prove some properties about certain monoids in a category, and to classify finite categories in terms of their endomorphism monoids. In particular, in this paper we study groups and monoids that contain a group. The classification of such monoids gives very interesting properties inside the category and gives us lots of additional information to conclude the number of categories with such monoids.
	
In this paper we present a construction theorem for grouplike categories, and we prove that every grouplike category should come from this construction.

\begin{definition} \label{unigen}
Let $L$ be a $(G_1^{*k_1},G_2^{*k_2})$-bimodule. We say that $L$ is \textit{$i$-unigen} if
$$
L \cdot f_i = e_i \cdot L := L_i 
$$
and $L_i \simeq G_1^{*i}$ as a left module and $L_i \simeq G_2^{*i}$ as a right module.
\end{definition}

\begin{lemma} \label{unigeniso}
Suppose $L$ is an $i$-unigen $(G_1^{*k_1},G_2^{*k_2})$-bimodule with $L_i$ as above, then there exists an isomorphism $G_1 \simeq G_2$ such that $L_i = G_1^{*i}$ as $(G_1^{*k_1},G_2^{*k_2})$-bimodule. If $i > 0$ the isomorphism is unique. If $i=0$ the isomorphism is well  defined up to inner automorphisms.
\end{lemma}

Because of this lemma we can assume that $G_1 = G_2 = G$. In this case we say that a $(G^{*k_1},G^{*k_2})$-bimodule $L$ is \textit{strongly $i$-unigen} if it is $i$-unigen and the isomorphism of Lemma \ref{unigeniso} can be taken as the identity. 

\begin{definition}
Let 
$$
M = \left( \begin{matrix}
	G^{*k_1} & L \\
	R & G^{*k_2}
\end{matrix} \right)
$$
be the algebraic matrix (Definition \ref{structuredmatrix}) of a grouplike category $\textsf{C}$, such that $I_1 = \{e_1,\hdots, e_{k_1}\}$ and $I_2 = \{f_1,\hdots, f_{k_2}\}$. We define $i^{max}(\textsf{C})$ to be the maximum index $i$ such that there exist $x,y$ such that $x\cdot y = e_i$ and $y \cdot x = f_i$.

\end{definition}
	Let $\textsf{C}$ be a category with one object $\{\alpha \}$ and endomorphism set $G^{*i}$.
		Let $x,y$ be two elements, then there exists a function
		$$
		f : \{x,y\} \rightarrow \alpha
		$$
		and we obtain a new category $\textsf{C}'$ such that $\textsf{C}' = f^{\alpha}(\textsf{C})$, $Ob(\textsf{C}') = \{x,y\}$ and $\textsf{C}'(x,y) = \textsf{C}(f(x),f(y))$.

		If we apply this construction to a grouplike monoid $G^{*i}$, we get a category associated to the matrix
		\begin{equation*} 
			\left( \begin{array}{cc}
				G^{*i} & G^{*i} \\
				G^{*i} & G^{*i} 
			\end{array} \right) \tag{$M^{(1)}$}
		\end{equation*}
	of similar copies of $G^{*i}$, such a category will be called \textit{groupoid-like}.
	\begin{theorem} \label{grouplikeconstruction}
		Let $G$ be a group and $G^{*i}$ a grouplike monoid of the form $G \cup I$ such that $I = \{e_1,\hdots, e_i\}$. Let $M^{(1)}$ be a matrix of the form
		\begin{equation} \label{mat1} 
			\left( \begin{array}{cc}
				G^{*i} & G^{*i} \\
				G^{*i} & G^{*i} 
			\end{array} \right) \tag{$M^{(1)}$}
		\end{equation}
		of similar copies of $G^{*i}$, let $\textsf{C}^{(1)}$ be the groupoid-like category associated to this matrix. Then we can extend the endomorphism sets $G^{*i}$ and we obtain a category $\textsf{C}^{(2)}$ associated to the matrix
		\begin{equation} \label{mat2} 
			\left( \begin{matrix}
				G^{*k_1} & G^{*i} \\
				G^{*i} & G^{*k_2} 
			\end{matrix} \right)\tag{$M^{(2)}$}
		\end{equation}
		such that for all $x \in G^{*i}, y \in G^{*k_1}$, we have $x \cdot y = y \cdot x \in G^{*i}$. Same for $G^{*k_2}$. And there exists $\textsf{C}^{(2)}$ a finite category associated to $M^{(2)}$.
		
		Now let $k_1,k_2 \geq i$, suppose that we have the matrix
		\begin{equation} \label{mat3} 
			\left( \begin{matrix}
				G^{*k_1} & L \\
				R & G^{*k_2} 
			\end{matrix} \right)\tag{$M^{(3)}$}
		\end{equation}
		such that $L$ and $R$ are strongly $i$-unigen. Then
		\begin{center}
			(\ref{mat1}) $\subseteq$ (\ref{mat2}) $\subseteq$ (\ref{mat3}).
		\end{center}
		And (\ref{mat3}) is a matrix of a unique grouplike category, denote it by $\textsf{C}^{(3)}$, such that $i^{max}(\textsf{C}^{(3)}) = i$ and for all $x \in L, y \in R$, $x \cdot y = e_i \cdot x \cdot y \cdot e_i$.
		
		We have
		$$
		\textsf{C}^{(1)} \subseteq \textsf{C}^{(2)} \subseteq \textsf{C}^{(3)}.
		$$
	\end{theorem}

\begin{theorem} \label{goaltheorem}
Every grouplike category comes from the construction described in Theorem \ref{grouplikeconstruction}.
\end{theorem}	
	
	\section{Preliminaries}
	In this section we introduce grouplike categories. These are groupoids (groups in the case of monoids) in which we add extra elements, having some specific properties, to their set of morphisms. The goal is to study the structure of these categories in order to make the classification problem easier and clearer.
	
	\begin{definition}
A \textit{monoid} $A$ is a set equipped with a binary operation $ \cdot : A \times A \rightarrow A$ such that $\cdot$ is associative and there exists an identity element $e$ such that for every element $a \in A$, the equations $e \cdot a = a$ and $a \cdot e = a$ hold.
	\end{definition}

\begin{definition}
	A \textit{bimodule} is a set with actions on the left and the right of the respective monoids, such that the actions commute i.e. $(g\cdot x)\cdot h  = g\cdot (x \cdot h)$. It can be seen as a category such that one of the sets of morphisms is empty, it's called an \textit{upper triangulated category}.
\end{definition}
	
\begin{definition} [\cite{leinster2002generalized}, \cite{leinster1999generalized}, \cite{leinster2-bimod}, \cite{Aurelio}, \cite{koslowski}] \label{structuredmatrix}
Let $\textsf{C}$ be a finite category. Then we get a matrix where the diagonal entries are monoids and the off diagonals entries are bimodules.
This matrix is called \textit{algebraic matrix}.
\end{definition}	
	
	\begin{definition}
		A \textit{semicategory} is a category without identity morphisms.
	\end{definition}
	
	\begin{lemma}
		Let $\textsf{C}$ be a semicategory, let $u\in \textsf{C}(x,y)$ with $x,y\in Ob(\textsf{C})$, then we can add a morphism $u'$ to $\textsf{C}(x,y)$ such that $u'\neq u$ and $u'$ duplicates $u$ for all composition operations. Then we get a new semicategory $\textsf{C}'$ with $Ob(\textsf{C}') = Ob(\textsf{C})$ and $Hom(\textsf{C}') = Hom(\textsf{C}) \cup \{u'\}$.
		
		On the other hand, we can also add the missing identities in $Ob(\textsf{C})$ to obtain a category $\textsf{B}$.
	\end{lemma}
	The previous lemma shows that we can add morphisms consecutively to a category and obtain a new category (provided that we add identities too). In our case, we define a category whose objects have grouplike endomorphism monoids, we only add elements consecutively to the monoids. The elements are idempotents and identities to the elements of the monoids. This means that each time we add an identity element, the previous identity is no longer an identity.

	We need the notion of bimodules in order to make the classification and counting problem easier.
	
	\begin{definition}
		Let $\textsf{A}$ be a semigroup. Whether or not $\textsf{A}$ has a multiplicative identity element, we let $e$ be a fresh element and
		\begin{center}
			$\textsf{A}^* := \textsf{A} \cup \{e\}$.
		\end{center}
		Then $\textsf{A}^*$ is a semigroup if the multiplication of $\textsf{A}$ is extended by stipulating $xe = ex = x$ for all $x \in \textsf{A}^*$. More generally, if $\textsf{A}$ is a semigroup, we recursively define semigroups $\textsf{A}^*$ for $n\in \mathbb{N}$ by:
		\begin{center}
			$\textsf{A}^{*0} = \textsf{A}$ \\
			$\textsf{A}^{*(n+1)} = (\textsf{A}^{*n})^*$.
		\end{center}
		If $\textsf{A}$ is a group, we say that the semigroups $\textsf{A}^{*n}$, $n \in \mathbb{N}$, are \textit{grouplike}.
	\end{definition}
		
	\begin{definition}
		We say that a category is called a \textit{grouplike category} with groups $G_{i}$ if its endomorphism monoids are \textit{grouplike monoids} of the form $G_{i}^{*k_i};~ k_i \in \mathbb{N}$.
	\end{definition}
	
	\begin{definition}
		A \textit{band} $S$ is an idempotent semigroup, i.e. for all $x \in S$, $x^2 = x$.\\
		A \textit{semilattice} is a commutative band.
	\end{definition}
	
	\begin{remark}
		The set of idempotents along with the group identity form a semilattice.
		\begin{proof}
			Let $I = \{e_0,\hdots, e_k\}$ be the set of idempotents where $e_0$ is a group identity. $I$ is a semilattice:
			\begin{itemize}
				\item commutative: for $i \leq j$, $e_i\cdot e_j = e_j \cdot e_i = e_i$.
				\item idempotent: $e_i \cdot e_i = e_i$.	
			\end{itemize}
		\end{proof}
	\end{remark}
	
	\section{Group action and the orbits of the sets of morphisms}
	\setlength{\parindent}{1em}
	Let $\textsf{C}$ be a category with $n$ objects $X_1,\ldots,X_n$. For each object $X_i, X_j$ there exists a monoid $A_i = \textsf{C}(X_i,X_i)$ and a monoid $A_j = \textsf{C}(X_j,X_j)$ and two operations
	$$
	l : A_i \times \textsf{C}(X_i,X_j) \rightarrow \textsf{C}(X_i,X_j)
	$$
	and
	$$
	r : \textsf{C}(X_i,X_j) \times A_j \rightarrow \textsf{C}(X_i,X_j)
	$$
	which represent the left monoid action of $A_i$ and the right monoid action of $A_j$ on the set of morphisms from $X_i$ to $X_j$. Note here that these are the entries of the matrix defined in Definition \ref{structuredmatrix}.
	
	In our work here, we take monoids of the form $\textsf{A}^{*n}$, specifically grouplike, which means each monoid contains a subgroup. This subgroup does not act on the whole set of morphisms from $X_i$ to $X_j$, but it acts on a subset of the previous set (it's important to note here that when we say group action we mean that the identity condition holds). In the following, we introduce how the group action works, and what are exactly the subsets that the group acts on. We will be considering categories with two objects.
	
	\begin{notation} \label{statement}
		We denote by $M(G^{*k_1}_{1},G^{*k_2}_{2},L,R)$ the matrix
		\begin{center}
			$M = \left( \begin{array}{cc}
				G_{1}^{*k_1} & L \\
				R & G_{2}^{*k_2}
			\end{array} \right)$
		\end{center}
		of monoids and bimodules such that $G_{1}^{*k_1} = G_1 \cup I_1$ and $G_{2}^{*k_2} = G_2 \cup I_2$, where $G_1$ and $G_2$ are groups and $I_1 =\{e_1,\hdots,e_{k_1}\}$ and $I_2 = \{f_1,\hdots f_{k_2}\}$ such that the elements of $I_1$ and $I_2$ satisfy (\ref{idsets}). We denote by $e_0$ and $f_0$ the identities of the groups $G_1$ and $G_2$. Let $Cat(M(G^{*k_1}_{1},G^{*k_2}_{2},L,R))$ be the set of grouplike categories associated to $M$ whose objects are are $X,Y$ such that the monoids and the bimodules are not empty. If $\textsf{C} \in Cat(M(G^{*k_1}_{1},G^{*k_2}_{2},L,R))$ then $\textsf{C}(X,X) = G_{1}^{*k_1}$, $\textsf{C}(X,Y) = L$, $\textsf{C}(Y,X) = R$ and $\textsf{C}(Y,Y) = G_{2}^{*k_2}$.
	\end{notation}
\begin{remark} \label{beta}
	More precisely if $M$ is an algebraic matrix i.e. a matrix of monoids and bimodules, then $Cat(M)$ is the set of pairs $(\textsf{C}, \beta)$ where $\textsf{C}$ is a category and $\beta$ is an isomorphism between the algebraic matrix of $\textsf{C}$ and $M$. We usually don't include this notation of $\beta$ in our discussion.
\end{remark}
	
	\begin{remark}
	Denote by $n_i$ the order of the group $G_i$. When we write $G^{*k_1}$ and $G^{'*k_2}$, this means that the groups $G$ and $G'$ have the same order $n$.
	\end{remark}
	
	\setlength{\parindent}{0em}
	\begin{lemma} \label{groupaction}
		Let $\textsf{C} \in Cat(M(G^{*k_1}_{1},G^{*k_2}_{2},L,R))$ (Notation \ref{statement}) be a grouplike category. Then:
		\begin{enumerate}
			\item $G_{1}$ acts on $G_{1} \cdot L$ and $R \cdot G_{1}$.
			\item $G_{2}$ acts on $L\cdot G_{2}$ and $G_{2} \cdot R$.
		\end{enumerate}
		\begin{proof}
			Let
			\begin{center}
				$ l: G_{1} \times (G_{1} \cdot L) \rightarrow G_{1} \cdot L$.
			\end{center}
			\begin{itemize}
				\item For all $x \in L$ and $g \in G_{1}$ we have $e_0 \cdot (g \cdot x) = g \cdot x$.
				\item For all $g_1, g_2, g_3 \in G_{1}$ and $x\in L$ we have $(g_2 \cdot g_3) \cdot (g_1 \cdot x) = g_2 \cdot (g_3 \cdot g_1 \cdot x)$.
			\end{itemize}
			Same for the other sets. \\
			Moreover, the orbit of the set $G_{1} \cdot L$ is itself. Indeed, the orbit of $G_{1} \cdot L$ is a subset $G_{1} \cdot (G_{1} \cdot L)$. It remains to prove the other direction. Let $g \cdot x \in G_{1} \cdot L$, then
			\begin{center}
				$g \cdot x = e_0 \cdot g \cdot x \in G_{1} \cdot (G_{1} \cdot L)$.
			\end{center}
		\end{proof}
	\end{lemma}
	
	Similarly to the case where we have groups as objects, we can conclude that the group action on these sets is free.
	
	\begin{prop}{\label{free_action}}
		Let $\textsf{C} \in Cat(M(G^{*k_1}_{1},G^{*k_2}_{2},L,R))$ (Notation \ref{statement}) be a grouplike category. The actions of the group $G_{1}$ on $G_{1} \cdot L$ and $R \cdot G_{1}$ and the group $G_{2}$ on $L \cdot G_{2}$ and $G_{2} \cdot R$ are all free.
		\begin{proof}
			Let $g_1,g_2,g \in G_{1}$ and $x \in L$. Suppose $g_1 \cdot ( g \cdot x) = g_2 \cdot (g \cdot x)$. If we multiply both sides by $y \in R$, we obtain:
			\begin{center}
				$g_1 \cdot g \cdot x \cdot y = g_2 \cdot g \cdot x \cdot y$
			\end{center}
			where $x\cdot y$ has 2 possibilities:
			\begin{enumerate}
				\item If $x \cdot y = e_i$ where $e_i \notin G_{1}$, then it's sufficient to multiply by $g^{-1}$ on both sides to prove that $g_1 = g_2$.
				\item If $x \cdot y = g_3$ where $g_3 \in G_{1}$, then if we multiply by $g_3^{-1}\cdot g^{-1}$ on both sides, we get $g_1 = g_2$.
			\end{enumerate}
		\end{proof}
	\end{prop}
	
	\begin{cor}\label{mult}
		Let $\textsf{C} \in Cat(M(G^{*k_1}_{1},G^{*k_2}_{2},L,R))$ (Notation \ref{statement}) be a grouplike category. The cardinal of each orbit in the set of morphisms is equal to the order of the group acting.
		\begin{proof}
			Use Proposition \ref{free_action}.
		\end{proof}
	\end{cor}
	
	This corollary could greatly help the enumeration problem here. Now that we know the number of possibilities in some blocks inside a category, then we can compute how many times the multiplication of morphisms appear to obtain non isomorphic copies of blocks.
	
	\begin{lemma}{\label{identity}}
		Let $\textsf{C} \in Cat(M(G^{*k_1}_{1},G^{*k_2}_{2},L,R))$ (Notation \ref{statement}) be a grouplike category. For all $x : X \rightarrow Y$ there exists at least one $y : Y\rightarrow X$ such that $x\cdot y = id_{G_{1}}$ and vice versa.
		\begin{proof}
			Let $x\in \textsf{C}(X,Y)$ and $y \in \textsf{C}(Y,X)$. Suppose $x \cdot y \neq id_{G_{1}}$:
			\begin{enumerate}
				\item If $x \cdot y = g \in G_{1}$, then
				\begin{center}
					$x \cdot y \cdot g^{-1} = id_{G_{1}}$.
				\end{center}
				\item If $x\cdot y$ is equal to an idempotent $e$ in $\textsf{C}(X,X)$ then
				\begin{eqnarray*}
					x \cdot y &= & e \\
					x \cdot y \cdot g &= & g~ ~;~~ g \in G_{1} \\
					x \cdot y \cdot g \cdot g^{-1} &= & id_{G_{1}}.
				\end{eqnarray*}
			\end{enumerate}
		\end{proof}
	\end{lemma}
	
	Now since we have two group actions on the sets of morphisms, we want to understand the relation between these actions over the morphisms sets.
	\begin{lemma}
		Let $\textsf{C} \in Cat(M(G^{*k_1}_{1},G^{*k_2}_{2},L,R))$ (Notation \ref{statement}) be a grouplike category. Then
		$$
		G_{1} \cdot x \cdot G_{2} \subseteq G_{1} \cdot x
		$$
		and 
		$$
		G_{1} \cdot x \cdot G_{2} \subseteq x \cdot G_{2}
		$$
		for all $x \in \textsf{C}(X,Y)$.
		\begin{proof}
			Let $g \cdot x \cdot h \in G_{1} \cdot x \cdot G_{2}$,
			\begin{eqnarray*}
				g \cdot x \cdot h &= & g \cdot x \cdot h \cdot e_2 ~~;~~ e_2 = id_{G_2}\\
				&= & g \cdot x \cdot h \cdot y \cdot x ~~;~~ y \in \textsf{C}(Y,X) ~~(\textrm{Lemma ~\ref{identity}})\\
				&\in & G_1 \cdot x.
			\end{eqnarray*}
			Same for the second inequality.
			
		\end{proof}
	\end{lemma}
	
	\begin{lemma}\label{isol}
		Let $\textsf{C} \in Cat(M(G^{*k_1}_{1},G^{*k_2}_{2},L,R))$ (Notation \ref{statement}) be a grouplike category. Let $x\in L$, the orbit of $x$ by the action of $G_{1}$ is the same orbit of $x$ by the action of $G_{2}$, i.e. $G_{1}\cdot x = x \cdot G_{2}$.
		\begin{proof}
			Suppose $n_1 \leq n_2$. From Lemma \ref{groupaction}, there is a group action by $G_{1}$ on $G_{1} \cdot L$, and by Corollary \ref{mult}
			\begin{center}
				$|G_{1} \cdot x| = n_1$.
			\end{center}
			Similarly, there is a group action by $G_{2}$ on $G_{1} \cdot L \cdot G_{2} \subseteq G_{1} \cdot L$, then again by Corollary \ref{mult} we have
			$$
			|g_1 \cdot x \cdot G_2| = n_2
			$$
			but
			$$
			g_1 \cdot x \cdot G_2 \subseteq G_1 \cdot x \cdot G_2 \subseteq G_1 \cdot x.
			$$
			
			Therefore, $n_2 \leq n_1$, but we have $n_1 \leq n_2$, then we obtain that $n_1$ should be equal to $n_2$ and $G_{1} \cdot x \cdot G_{2} = G_{1} \cdot x$ and $G_{1} \cdot x \cdot G_{2} = x\cdot G_{2}$, hence
			\begin{center}
				$G_{1} \cdot x = x \cdot G_{2}$ for all $x \in L$.
			\end{center}
		\end{proof}
	\end{lemma}

	\begin{prop}
		Let $\textsf{C} \in Cat(M(G^{*k_1}_{1},G^{*k_2}_{2},L,R))$ (Notation \ref{statement}) be a grouplike category. For all $x, y \in L$ we have
		\begin{center}
			$G_{1} \cdot x = G_{1} \cdot y$.
		\end{center}
		\begin{proof}
			Let $g \cdot x \in G_{1} \cdot x$,
			\begin{eqnarray*}
				g \cdot x &= & (g \cdot e_1) \cdot x ~ ~; ~~ e_1 = id_{G_{1}} \\
				&= & g \cdot (e_1 \cdot x) \\
				&= & g \cdot e_1 \cdot x \cdot e_2~~ (Lemma ~\ref{isol}) \\
				&= & g \cdot e_1 \cdot x \cdot z \cdot y \in G_{1} \cdot y~~ (Lemma ~\ref{identity})
			\end{eqnarray*}
			and since $|G_{1} \cdot x| = |G_{1} \cdot y|$, hence equality.
		\end{proof}
	\end{prop}
	
	\begin{prop} \label{compoforbits}
		Let $\textsf{C} \in Cat(M(G^{*k_1},G^{'*k_2},L,R))$ (Notation \ref{statement}) be a grouplike category. Then the multiplication of the elements in the orbit of $L$ and the elements in the orbit of $R$ is the group $G$, i.e.
		$$
		G \cdot \textsf{C}(X,Y) \cdot \textsf{C}(Y,X) \cdot G = G
		$$
		and
		$$
		G' \cdot	\textsf{C}(Y,X) \cdot \textsf{C}(X,Y) \cdot G' = G'.
		$$
		
		\begin{proof}
			\begin{itemize}
				\item $\subseteq$ : evident.
				\item $\supseteq$ : Let $g \in G$
				\begin{center}
					$ g = g \cdot e_0 \cdot e_0 = g \cdot x \cdot y \cdot e_0	$ (Lemma \ref{identity})
				\end{center}
				where $x : X \rightarrow Y$ and $y : Y \rightarrow X$.
			\end{itemize}
		\end{proof}
	\end{prop}

	\begin{definition} \label{isolated}
		Let $\textsf{C}$ be a category with objects $\{X,Y\}$. Suppose that $g\cdot f$ is never equal to an identity for all $f,g \in Hom(\textsf{C})$, then we can always reduce the category $\textsf{C}$ to a new semi-category $\textsf{C}'$ such that $Ob(\textsf{C}') = Ob(\textsf{C})$ and $Hom(\textsf{C}') = Hom(\textsf{C}) \setminus \{id_\textsf{C}\}$.
	\end{definition}
	Similarly, we can eliminate morphisms other than identities and still obtain a new semi-category. 
	%
	
	The reason we want to eliminate some morphisms is because when we take a category whose objects have grouplike endomorphism monoids, then we could restrict the category to a smaller category with only groups as objects and the orbit of the set of morphisms. Following this technique leads to proving some matrix properties about the coefficients on the off-diagonals. It will also clarify how such categories are built.
	
	From the above lemmas and propositions, we conclude that we can divide each set of morphisms into 2 sets: the orbit of the set and the other elements that are not inside the orbit. 
	
	\subsection{Groupoids}
	\begin{definition}
		In category theory, a \textit{groupoid} generalizes the notion of group in several equivalent ways. A groupoid can be seen as a:
		\begin{itemize}
			\item group with a partial function replacing the binary operation;
			\item category in which every morphism is invertible. A groupoid with only one object is a usual group.	
		\end{itemize}
	\end{definition}
	
	In the previous sections, we have seen the structure of a grouplike category with two objects. From this these categories, we can extract a subcategory, which is exactly a groupoid.
	
	\begin{theorem}
		In every grouplike category with two objects, there is a sub-semicategory that is a groupoid with two objects, whose groups are the groups 
		of the grouplike monoids.
		\begin{proof}
			Let $\textsf{C}$ be a grouplike category with objects $X,Y$. Let $\textsc{C}(X,X) = G^{*k_1}, \textsc{C}(Y,Y) = G^{'*k_2}, \textsc{C}(X,Y) = L$ and $\textsc{C}(Y,X) = R$. $\textsf{G}$ the sub-semicategory of $\textsf{C}$ is of the form:
			\begin{itemize}
				\item $Ob(\textsf{G}) = \{G, G^{'}\}$;
				\item $Mor(\textsf{G}) = \{G, G^{'}, G \cdot L = L \cdot G^{'}, G^{'} \cdot R = R \cdot G \}$;
				\item Identities of $\textsf{G}$ are $1_{G}$ and $1_{G^{'}}$;
				\item Let $x,y \in \textsf{G}$, $x \circ y = x \cdot y$ with $\circ$ associative.
			\end{itemize}
		\end{proof}
	\end{theorem}
	
	\begin{cor}
		Let $\textsf{C} \in Cat(M(G^{*k_1},G^{'*k_2},L,R))$ (Notation \ref{statement}) be a grouplike category. The category $\textsf{C}$ determines an isomorphism between $G$ and $G'$. The isomorphism is well  defined up to inner automorphisms.
		\begin{proof}
			Evident.
		\end{proof}
	\end{cor}
	
	\subsection{The sets of idempotents}
	In this section we study the role of the idempotent elements in the monoids, the idea is to interpret their action on the sets of morphisms. 
	
	Let $\textsf{C} \in Cat(M(G_1^{*k_1},G_2^{*k_2},L,R))$ (Notation \ref{statement}) be a grouplike category. Let:
		$$
		K_{1} = \{ e_0, e_1, \hdots , e_{k_1} \}
		$$
		and
		$$
		K_{2} = \{ f_0, f_1, \hdots , f_{k_2} \}
		$$
	be the sets of idempotents of $G_1^{*k_{1}}$ and $G_2^{*k_{2}}$ where $e_0 = 1_{G_1}$ and $f_0 = 1_{G_2}$. And let:
	\begin{center}
		$ L_{ij} := \{x \in L \mid e_i \cdot x = x ~\textrm{and}~ x \cdot f_j = x \} = e_i \cdot L \cdot f_j$\\
		and \\
		$ R_{ji} := \{y \in R \mid y\cdot e_i = y ~\textrm{and}~ f_j \cdot y = y \} = f_j \cdot R \cdot e_i$
	\end{center}
	be the sets of morphisms that are fixed by $e_i$ and $f_j$.
	
	\begin{notation}
		$L_i : = L_{ii}$ and $R_j : = R_{jj}$.
	\end{notation} 
	
	In general we have
	\begin{lemma} \label{intersection}
		$L_{ij} = e_i \cdot L \cap L \cdot f_j$ and $R_{ji} = f_j \cdot R \cap R \cdot e_i$.
		\begin{proof}
			We always have
			$$
			L_{ij} \subseteq e_i \cdot L ~ \textrm{and} ~ L_{ij} \subseteq L \cdot f_j
			$$
			then
			$$
			L_{ij} \subseteq e_i \cdot L \cap L \cdot f_j.
			$$
			Now let $x \in e_i \cdot L \cap L \cdot f_j$ then $e_i \cdot x = x$ and $x \cdot f_j = x$ hence $x \in L_{ij}$.
		\end{proof}
	\end{lemma}
	
	In the following lemma, we present some properties of the multiplication of the sets of morphisms by idempotent elements.
	
	\begin{lemma}\label{idemprops}
		Let $\textsf{C} \in Cat(M(G_1^{*k_1},G_2^{*k_2},L,R))$ (Notation \ref{statement}) be a grouplike category, we will denote by $x,x'$ elements of $L$ and by $y,y'$ elements of $R$.
		\begin{enumerate}
			\item $e_i \cdot x = x$ and $x \cdot f_j = x$ for all $e_i \in K_1$, $f_j \in K_2$ and $x \in L_0$.
			\item $e_j \cdot (e_i \cdot x) = e_i \cdot x$ and $(x \cdot f_i) \cdot f_j = x \cdot f_i$ for all $x\in L$ and $i \leq j$.
			\item $x \cdot y \in G_1$ and $y \cdot x \in G_2$ for all $x \in L_0$ or $y \in R_0$.
			\item If $x\cdot y \in K_1 \setminus \{e_0\}$ then $y \cdot x \in K_2\setminus \{f_0\}$ and vice versa.\label{condidemp}
			\item If $x \cdot y = e_0$ then $y \cdot x = f_0$. (This result is also proved in Lemma \ref{maxij}). \label{condidempzero}
			\item \label{cond7} If $x\cdot y = e_i$ and $y\cdot x = f_j$ then $x \cdot f_j = e_i \cdot x$ and $y \cdot e_i = f_j \cdot y$.
			
			In this case we can always assume that $x \cdot f_j = e_i \cdot x = x$ and $y \cdot e_i = f_j \cdot y = y$ (because $x$ and $y$ could be replaced with $e_i \cdot x = x \cdot f_j$ and $ y \cdot e_i=f_j \cdot y$ respectively). \label{lemma6}
		\end{enumerate}
		\begin{proof}
			\begin{enumerate}
				\item Let $x \in L_0$ and $e_i \in K_1$, we have:
				\begin{center}
					$e_i \cdot x = e_i \cdot (e_0 \cdot x) = (e_i \cdot e_0) \cdot x = e_0 \cdot x = x$.
				\end{center}
				Same for $x \cdot f_j$.
				\item $e_i \cdot e_j = e_j \cdot e_i = e_i$ for all $i \leq j$ then $e_j \cdot (e_i \cdot x) = x$.
				\item Let $g\in G_1$ such that $g \cdot x = g\cdot x'$ then $g^{-1}\cdot g \cdot x = g^{-1}\cdot g\cdot x'$ then $g'\cdot x = g' \cdot x'$.
				\item Suppose that $x \cdot y = e_i$ then $x \cdot y \cdot e_0 = e_0$ then $x \cdot y = e_0$ (because $y \cdot e_0 = y$) contradiction.
				\item Let $x \in L$ and $y\in R$ such that $x \cdot y = e_i \in K_1\setminus\{e_0\}$. Suppose that $y \cdot x = g\in G_2$. Then
				\begin{eqnarray*}
					(x\cdot y) \cdot (x \cdot y) &=& e_i \\
					x \cdot g \cdot y &=& e_i
				\end{eqnarray*}
				but $x \cdot g \cdot y \in G_1$ and $e_i \in K_1\setminus\{e_0\}$, it means that $e_i \in G_1$, hence contradiction. Therefore, $y \cdot x \in K_2 \setminus\{f_0\}$.
				\item From part (\ref{condidemp}), we can see that if $x\cdot y \in G_1$ then $y \cdot x \in G_2$. Suppose that $x \cdot y = e_0$ and $y \cdot x = g \in G_2$. Proving that $g = f_0$. We have
				\begin{eqnarray*}
					y \cdot x &=& g \\
					\Rightarrow x \cdot y \cdot x &=& x \cdot g \\
					\Rightarrow y \cdot e_0 &=& g \cdot y \\
					\Rightarrow f_0 \cdot (y \cdot e_0) &=& g \cdot (y \cdot e_0)~~(\textrm{because}~ y\cdot e_0~\textrm{and}~ g\cdot y~\textrm{are in}~G_2 \cdot R).
				\end{eqnarray*}
				By the free action of $G_2$ on $G_2 \cdot R$, we obtain that $g = f_0$.
				\item $x\cdot y \cdot x = e_i\cdot x$ thus $x \cdot f_j = e_i \cdot x$\\
				$y \cdot x \cdot y = f_j\cdot y $ thus $y \cdot e_i = f_j \cdot y$.
				
				In addition, by part (b) $(x \cdot f_j )\cdot f_j = x \cdot f_j$ and $e_i \cdot (e_i \cdot x ) = e_i \cdot x$, then we can assume that $x \cdot f_j = e_i \cdot x = x$.
			\end{enumerate}
		\end{proof}
	\end{lemma}

From Lemma \ref{idemprops} (\ref{condidemp}) (\ref{condidempzero}), we see that if one of the multiplications is an idempotent then the other way around should be an idempotent as well. That's why in the following, we study the structure of the category whenever we have two elements such that their multiplications are idempotents.
	
	\begin{lemma} \label{maxij}
		Let $\textsf{C} \in Cat(M(G_1^{*k_1},G_2^{*k_2},L,R))$ (Notation \ref{statement}) be a grouplike category. Let $i$ be the maximum index such that there exist $x\in L$ and $y\in R$; $x \cdot y = e_i$. And let $j'$ be the maximum index such that there exist $x'\in L$ and $y'\in R$; $y' \cdot x' = f_{j'}$. We can assume following Lemma \ref{idemprops} (\ref{lemma6}) that $e_i \cdot x = x$, $y \cdot e_i = y$, $x' \cdot f_{j'} = x'$ and $f_{j'} \cdot y' = y'$.
		\begin{enumerate}
			\item 	If $i=0$ then $j' =0$ and vice versa.
			\item If $i,j' \geq 1$ then $x\cdot y' = e_i$ and $y' \cdot x = f_{j'}$. In addition, $e_i \cdot x = x \cdot f_{j'} = x$ and $y' \cdot e_i = f_{j'} \cdot y' = y'$.
		\end{enumerate}
		\begin{proof}
			For part (1), suppose $i =0$, in this case by Lemma \ref{idemprops} (\ref{condidemp}) we have $x' \cdot y'$ is an idempotent, and by maximality of $i$ it has to be $e_0$. Therefore
			$$
			f_{j'} = f_{j'}^2 = y' \cdot x' \cdot y' \cdot x' 
					 = y' \cdot e_0 \cdot x' \in G_1.
			$$
			Therefore, $j' = 0$.
			
			For part (2), suppose $i',j' \geq 1$, then from Lemma \ref{idemprops} (\ref{condidemp}), suppose that $y \cdot x = f_j $ and $x' \cdot y' = e_{i'}$. From our assumption following (Lemma \ref{idemprops} (\ref{lemma6})), we get
			$$
			f_j \cdot y = y \cdot e_i = y~ \textrm{and} ~x \cdot f_j = e_i \cdot x = x
			$$
			and
			$$
			f_{j'} \cdot y' = y' \cdot e_{i'} = y' ~ \textrm{and} ~ x' \cdot f_{j'} = e_{i'} \cdot x' = x'.
			$$
			By maximality we have $i \geq i'$ and $j' \geq j$. We have
			\begin{center}
				$x \cdot f_{j'} = (x\cdot f_j) \cdot f_{j'} = x \cdot (f_j \cdot f_{j'}) = x \cdot f_j = x$.
			\end{center}
			We obtain that
			\begin{center}
				$x = x \cdot f_j = x \cdot f_{j'} = x \cdot (y' \cdot x')$.
			\end{center}
			Then
			\begin{eqnarray*}
				e_i = x\cdot y = (x \cdot y' \cdot x') \cdot y = (x\cdot y') \cdot (x' \cdot y).
			\end{eqnarray*}
		Therefore $x \cdot y' \notin G_1$ since $i \geq 1$, then there exists $m \leq i$ (by maximality of $i$) such that $x \cdot y' =e_m$. Then $e_i = e_m \cdot (x' \cdot y')$.
		
			Which implies that $e_m = e_m \cdot e_i = e^2_m \cdot (x' \cdot y) = e_m \cdot (x' \cdot y) = e_i$. Hence $e_m = e_i$ and $x \cdot y' = e_i$.
			
			Similarly we can prove that $y' \cdot x = f_{j'}$ and $y' \cdot e_i = y'$.
		\end{proof}
	\end{lemma}
	
	\begin{cor}
		From Lemma \ref{maxij}, we have $x = x'$ and $y = y'$.
		\begin{proof}
			$x = e_i \cdot x = x \cdot y' \cdot x = x \cdot f_{j'} \cdot x \cdot y' \cdot x' = e_i \cdot x'$.
			
			Where $x' = e_{i'} \cdot x' = e_{i'} \cdot e_i \cdot x' = e_i \cdot e_{i'} \cdot x' = e_i \cdot x'$. Then $x = x'$.
			
			Similarly we prove that $y = y'$.
		\end{proof}
	\end{cor}
	
	\begin{prop}
		Let $\textsf{C} \in Cat(M(G_1^{*k_1},G_2^{*k_2},L,R))$ (Notation \ref{statement}) be a grouplike category. Let $i$ and $j$ be the maximum elements such that there exist $x\in L$ and $y\in R$; $x \cdot y = e_i$ and $y \cdot x = f_j$. Then
		$$
		L_{ij} = e_i \cdot L = L \cdot f_j
		$$
		and
		$$
		R_{ji} = f_j \cdot R = R \cdot e_i.
		$$ 
		\begin{proof}
			We want to prove that $e_i \cdot L \cdot f_j = e_i \cdot L = L \cdot f_j$.
			
			We always have that $e_i \cdot L \cdot f_j \subseteq e_i \cdot L$ and $e_i \cdot L \cdot f_j \subseteq L \cdot f_j$. We prove the other direction.
			
			Let $x'\in L$
			\begin{eqnarray*}
				e_i x' &=& e_i e_i x'\\
				&=& e_i (x y) x' \\
				&=& e_i x (yx') \\
				&=& e_i x f_m \\
				&=& e_i \underbrace{x f_m}_{\in L} f_j ~\textrm{($f_mf_j = f_m$~ because $j$ is the maximum)}~ \\
				&\in& e_i\cdot L \cdot f_j.
			\end{eqnarray*}	
			Similarly we prove the others.	
		\end{proof}
	\end{prop}

	\begin{theorem} \label{theorem}
		Let	$\textsf{C} \in Cat(M(G_1^{*k_1},G_2^{*k_2},L,R))$ (Notation \ref{statement}) be a grouplike category.
		
		Suppose that $i$ and $j$ are the maximum elements (in the sense of Lemma \ref{maxij}) such that $x\cdot y = e_i$ and $y \cdot x = f_j$. By Lemma \ref{idemprops} (\ref{cond7}) we can assume that
		\begin{center}
			$f_j \cdot y = y \cdot e_i = y$ and $x \cdot f_j = e_i \cdot x = x$
		\end{center}
		Then we can construct a sub-semicategory $\textsf{C}'$ of the form
		\begin{center}
			$ \left(
			\begin{array}{cc}
				G_1^{*i} & L_{ij} \\
				R_{ji} & G_2^{*j}
			\end{array}
			\right)$
		\end{center}
		where $L_{ij} = e_i \cdot L \cdot f_j$ and $R_{ji} = f_j \cdot R \cdot e_i$, such that $\textsf{C}'$ is the maximum sub-semicategory of this form.
		\begin{proof}
			\begin{itemize}
				\item \textbf{$\textsf{C}'$ is a category:}
				\begin{itemize}
					\item Objects: $X$ and $Y$
					\item Morphisms: $G_1^{*i}, L_{ij}, R_{ji}$ and $G_2^{*j}$
					\item Composition:
					\begin{itemize}
						\item Let $z \in L_{ij}$ and $e' \in K_i$ then we can write $z = e_i \cdot \tilde{z} \cdot f_j$, then
						\begin{center}
							$e' \cdot z = e' \cdot (e_i \cdot \tilde{z} \cdot f_j) = (e' \cdot e_i) \cdot \tilde{z} \cdot f_j = e_i \cdot (e' \cdot \tilde{z})\cdot f_j \in L_{ij}$
						\end{center}
						\item Let $z\in L_{ij}$ and $w \in R_{ji}$ then $z = e_i \cdot \tilde{z} \cdot f_j$ and $w = f_j \cdot \tilde{w} \cdot e_i$, then
						\begin{eqnarray*}
							z \cdot w &=& (e_i \cdot \tilde{z} \cdot f_j) \cdot (f_j \cdot \tilde{w} \cdot e_i) \\
							&=& e_i \cdot (\tilde{z} \cdot f_j \cdot f_j \cdot \tilde{w}) \cdot e_i \\
							&=& e_i \cdot (\tilde{z} \cdot f_j \cdot \tilde{w}) \cdot e_i \in G_1^{*i}
						\end{eqnarray*}
					\end{itemize}
					\item Identities: $e_i$ and $f_j$
				\end{itemize}
				\item \textbf{$\textsf{C}'$ is the maximum category of this form:}\\
				If we have two elements $x'$ and $y'$ outside of $L_{ij}$ and $R_{ji}$ such that $x' \cdot y' = e_{i'}$ and $y' \cdot x' = f_{j'}$ then by Lemma \ref{maxij} we have $x \cdot y' = e_i$ and $y' \cdot x = f_{j'}$. This is a contradiction with the maximality of $i$ and $j$. Hence $\textsf{C}'$ is maximum.
			\end{itemize}
		\end{proof}
	\end{theorem}
	
\begin{cor} \label{remark}
	The monoids $G_1^{*i}$ and $G_2^{*j}$ are isomorphic. We call $\textsf{C}'$ a \textit{groupoid-like} category.
	\begin{proof}
		There exist $x,y$ such that $x \cdot y = e_i$ and $y \cdot x = f_j$ the identities. It follows that $i = j$ because the groups $G_1$ and $G_2$ are isomorphic.
	\end{proof}
\end{cor}
	
	\begin{prop} \label{prop}
	Let
			\begin{center}
		$M = \left( \begin{array}{cc}
			G_1^{*k_1} & L \\
			R & G_2^{*k_2}
		\end{array} \right)$
	\end{center}
be a matrix of a grouplike category $\textsf{C}$. By Theorem \ref{theorem}, we can construct a sub-semicategory $\textsf{C}'$ associated to the matrix
		\begin{center}
			$M' = \left( \begin{array}{cc}
				G_1^{*i} & L_{i} \\
				R_{i} & G_2^{*i}
			\end{array} \right)$
		\end{center}
		of monoids and bimodules. For all $x \in L$ and all $y \in R$ we have 
		$$
		x \cdot y = (e_i \cdot x) \cdot (y \cdot e_i) \in G_1^{*i}~ \textrm{and}~ y\cdot x = (f_i \cdot y) \cdot (x \cdot f_i) \in G_2^{*i}.
		$$.
		\begin{proof}
			Suppose that $x \cdot y \in G_1^{*k_1} \setminus G_1^{*i}$ then there exists $e_k \notin G_1^{*i}$ such that $x \cdot y = e_k$. This is a contradiction with the maximality of $i$.
		\end{proof}
	\end{prop}
	
	\textbf{Conclusion:} Let
	$$
	M= \left( \begin{array}{cc}
		G_{1}^{*k_1} & L \\
		R & G_{2}^{*k_2}
	\end{array} \right)
	$$
	be a matrix of a grouplike category $\textsf{C}$ whose objects $X$ and $Y$. Then
	\begin{enumerate}
		\item $G_{1} \simeq G_{2}$.
		\item There exists $i^{max} = max(i)$ such that there exist $x,y$ and $e_{i^{max}} = x\cdot y$.
		
		There exists $j^{max} = max(j)$ such that there exist $x,y$ and $f_{j^{max}} = y\cdot x$.
		
		\item $i^{max} = j^{max}$, and we obtain the following matrix
		$$
		\left( \begin{array}{cc}
			G_1^{*i^{max}} & L_{i^{max}j^{max}} \\
			R_{j^{max}i^{max}} & G_2^{*j^{max}}
		\end{array} \right)
		$$
		of a sub-semicategory $\textsf{C}' $ of $\textsf{C}$ whose objects $X$ and $Y$ are isomorphic. Thus
		$$
		G_1^{*i^{max}} \simeq G_2^{*j^{max}} ~\textrm{as monoids}
		$$
		and
		$$
		G_1^{*i^{max}} \simeq L_{i^{max}j^{max}} \simeq R_{j^{max}i^{max}} ~\textrm{as bimodules}.
		$$
		\item \label{item} The bimodule $L$ has the property
		$$
		L_{i} = L \cdot f_i = e_i \cdot L \simeq G_1^{*i} \simeq G_2^{*i}
		$$
		i.e. $L$ is $i$-unigen.
		
		and the bimodule $R$ has the property
		$$
		R_i = f_i \cdot R = R \cdot e_i \simeq G_2^{*i} \simeq G_1^{*i}
		$$
		i.e. $R$ is $i$-unigen. Where $i = i^{max} = j^{max}$.
		\item For all $x \in L$, $e_i \cdot x = x \cdot f_i$.
		\item The isomorphisms in \ref{item} are inverses.
		\item The multiplications of the elements of $L$ by the elements of $R$ are determined once we fix $x$ and $y$.
	\end{enumerate}
	
	\begin{theorem}\label{grouplikecounting}
		Let $G_1, G_2$ be two groups, $k_1,k_2$, $L$, $R$ be $(G_1^{*k_1},G_2^{*k_2})$-bimodules, and $i$; $i \leq k_1$, $i\leq k_2$ such that
		\begin{equation}\label{P} \tag{P(i)}
			L_i = L\cdot f_i = e_i \cdot L ~\textrm{and}~
			R_i = R \cdot e_i = f_i \cdot R
		\end{equation}
		and such that there is $x_i\in L_i$ such that $x_i$ determines an isomorphism
		$$
		\begin{array}{ccc}
			G_1^{*i} & \rightarrow & L_i \\
			g & \mapsto & gx_i
		\end{array}
		$$
		and an isomorphism
		$$
		\begin{array}{ccc}
			G_2^{*i} & \rightarrow & L_i \\
			h & \mapsto & x_ih
		\end{array}.
		$$	
		Similarly for $y_i \in R_i$, determines an isomorphism $G_2^{*i} \simeq R_i \simeq G_1^{*i}$.
		
		The isomorphisms $G_1^{*i} \simeq L_i \simeq G_2^{*i}$ and $G_2^{*i} \simeq R_i \simeq G_1^{*i}$ are assumed to be inverses for \ref{P}.
		
		Then we get a category $\textsf{C}$ with algebraic matrix
		$$
		\left( \begin{array}{cc}
			G_1^{*k_1} & L \\
			R & G_2^{*k_2}
		\end{array} \right)
		$$
		such that $i = i^{max} = j^{max}$, $x_i \cdot y_i = e_i$ and $y_i \cdot x_i = f_i$ ($L$ and $R$ are $i$-unigen).
		
		For $i \geq 1$, then the choice of $x_i$ is unique and hence the category is unique.
		
		For $i = 0$, then the choice of $x_0$ is not unique
		but the category is unique once $x_0$ and $y_0$ are fixed.
		\begin{proof} The multiplications $G_1^{*k_1} \times L \rightarrow L$ and $ L \times G_2^{*k_2} \rightarrow L$ are given by the bimodule structure of $L$. Similarly for $R$.
			
			Let $x' \in L$, $y' \in R$,
			$$
			x'y' := \underbrace{(e_ix')}_{\in L_i}\underbrace{(y'e_i)}_{\in R_i} \in G_1^{*i}.
			$$
			Similarly for $y'x'$. We can check that the multiplication is associative.
			
			If we fix the matrix
			$$
			M = \left( \begin{array}{cc}
				G_1^{*k_1} & L \\
				R & G_2^{*k_2}
			\end{array} \right)
			$$
			and we take $i^{\textrm{biggest}} = max\{i \mid \textrm{\ref{P} holds}\}$, then for the cases $i \geq 1$
			$$
			Cat_{i^{max} \geq 1}(M) = \{1 \leq i \leq i^{\textrm{biggest}}\}.
			$$
			
For $i =0$, the choices of $x_0$ and $y_0$ are not unique, see Remark \ref{center}.
		\end{proof}
	\end{theorem}

\begin{remark}
The condition \ref{P} is what we call $i$-unigen in Definition \ref{unigen} (Thanks to Carlos Simpson for suggesting the terminology $i$-unigen).
\end{remark}

\begin{remark} \label{center}
When $i=0$, once we fix $x_0$, there are maybe several choices for $y_0$ that lead to inverse isomorphisms. The set of choices of $y_0$ is given by the center of the group. The cardinal of the set of pairs $(\textsf{C}, \beta)$ in $Cat(M(G_1^{*k_1},G_2^{*k_2},L,R))$ with $i^{max} = 0$ is equal to the cardinal of the center of the group.
\end{remark}

\noindent
{\em Proof of Lemma \ref{unigeniso}}. Let $L$ be an $i$-unigen bimodule. Let $x \in L_i$ such that $G_1^{*i} \cdot x = L_i$, where $G_1^{*i} \cdot x \simeq G_1^{*i}$. Then $x \cdot G_2^{*i} = L_i$:

Suppose $x \cdot g = x\cdot g'$, $g \neq g'$ in $G_2^{*i}$. Then since $G_1^{*i} \cdot x = L_i$, we get $y \cdot g = y \cdot g'$ for all $y \in L_i$, then $g = g'$. Then
$$
\{x\} \times G_2^{*i} \hookrightarrow L_i
$$
but they have the same cardinal, then
$$
\{x\} \times G_2^{*i} \simeq L_i.
$$
We define an isomorphism
$$
\phi : G_1^{*i} \rightarrow G_2^{*i}
$$
as follows

If $g \in G_1^{*i}$ and  $g \cdot x \in L_i$ then $g\cdot x = x \cdot h$; $h \in G_2^{*i}$. Set
\begin{center}
	$\phi(g) := h$ and $g \cdot x = x \cdot \phi(g)$.
\end{center}
Then
\begin{eqnarray*}
	(gg') \cdot x &=& g \cdot (g' \cdot x) \\
	&=& g \cdot (x \cdot \phi(g')) \\
	&=& (g \cdot x) \cdot \phi(g')\\
	&=& (x \cdot \phi(g)) \cdot \phi(g') \\
	&=& x \cdot (\phi(g) \phi(g')) \\
	&=& x \cdot (\phi(gg')).
\end{eqnarray*}
By uniqueness of $h$, we have $\phi(gg') = \phi(g) \phi(g')$, and Lemma \ref{unigeniso} is proved.\\
$~~~~~~~~~~~~~~~~~~~~~~~~~~~~~~~~~~~~~~~~~~~~~~~~~~~~~~~~~~~~~~~~~~~~~~~~~~~~~~~~~~~~~~~~~~~~~~~~~~~~~~~~~~~~~~~~~~~~~~~~~~~~~~$\qedsymbol{}\\

\noindent
{\em Proof of Theorem \ref{grouplikeconstruction}}. Using Theorem \ref{grouplikecounting} we prove Theorem \ref{grouplikeconstruction}, as we have the same construction of a category. In Theorem \ref{grouplikecounting} we start by choosing $i^{max}=j^{max}$ to get to the algebraic matrix (\ref{mat3}). 

For the multiplication table of $\textsf{C}^{(3)}$, the multiplication of $G^{*k_1}$ on $L$ and $R$ is given by the bimodule structure. It remains to find the maps
\begin{center}
	$L \times R \rightarrow G^{*k_1}$ and $R \times L \rightarrow G^{*k_2}$.
\end{center}
Let $x \in L$ and $y\in R$,
$$
x \cdot y = e_i \cdot x \cdot y \cdot e_i
$$
where $e_i \cdot x$ and $y \cdot e_i$ are in $L_i$ and $R_i$ and these compositions are given by $\textsf{C}^{(1)}$. And as in the conclusion part (7), we get the uniqueness of the category.\\
$~~~~~~~~~~~~~~~~~~~~~~~~~~~~~~~~~~~~~~~~~~~~~~~~~~~~~~~~~~~~~~~~~~~~~~~~~~~~~~~~~~~~~~~~~~~~~~~~~~~~~~~~~~~~~~~~~~~~~~~~~~~~~~$\qedsymbol{}\\

\noindent
{\em Proof of Theorem \ref{goaltheorem}}. From Theorem \ref{theorem} and Proposition \ref{prop} we can prove Theorem \ref{goaltheorem}. As in Remark \ref{remark} we prove that $i$ and $j$ should be the same, then the algebraic matrix obtained is
		\begin{center}
	$M = \left( \begin{array}{cc}
		G^{*i} & L_{i} \\
		R_{i} & G^{*i}
	\end{array} \right)$
\end{center}
where $L_i$ and $R_i$ are isomorphic to $G^{*i}$ (conclusion part (3)). Then $M$ is the matrix (\ref{mat1}).

We conclude that if we have a grouplike category then the two bimodules $L$ and $R$ are $i$-unigen and the resulting isomorphisms between these two bimodules are inverses. If we then identify the groups via these isomorphisms, we can say that $L$ and $R$ become strictly $i$-unigen. Then we get the structure described in Theorem \ref{grouplikeconstruction} and Theorem \ref{goaltheorem} is proved.\\
$~~~~~~~~~~~~~~~~~~~~~~~~~~~~~~~~~~~~~~~~~~~~~~~~~~~~~~~~~~~~~~~~~~~~~~~~~~~~~~~~~~~~~~~~~~~~~~~~~~~~~~~~~~~~~~~~~~~~~~~~~~~~~~$\qedsymbol{}\\
	
	\begin{notation}
Let 
$$
N^{(1)} = \left( \begin{matrix}
	A & L \\
	\emptyset & B
\end{matrix} \right)
$$
be a matrix of bimodule $L$, and let $\textsf{C}^{(1)}$ be a category associated to $N^{(1)}$.
 
Similarly let
$$
N^{(2)} = \left( \begin{matrix}
	A & \emptyset \\
	R & B
\end{matrix} \right)
$$
be a matrix of bimodule $R$, and let $\textsf{C}^{(2)}$ be a category associated to $N^{(2)}$.

We denote by $N = N^{(1)} \cup N^{(2)}$ the matrix of the form
$$
\left( \begin{matrix}
	A & L \\
	R & B
\end{matrix} \right)
$$
of monoids and bimodules $L$ and $R$ in $\textsf{C}^{(1)}$ and $\textsf{C}^{(2)}$. We denote by $\textsf{C}$ a category associated to $N$.
	\end{notation}
	
	\begin{remark}
		The condition that the isomorphisms should be inverses in Theorem \ref{grouplikecounting} is very important to obtain the grouplike category. For example, consider the matrix
		$$
		N^{(1)} = \left( \begin{matrix}
			\mathbb{Z}_3 & L \\
			\emptyset & \mathbb{Z}_3
		\end{matrix} \right)
		$$
		where
		$$
		\begin{matrix}
			\mathbb{Z}_3 = \{1,2,3\} & L = \{4,5,6\} & \mathbb{Z}_3 = \{7,8,9\}	
		\end{matrix} 
		$$
		such that
		$$
		\mathbb{Z}_3 \cdot L = L \cdot \mathbb{Z}_3 = \left( \begin{matrix}
			4 & 5 & 6 \\
			5 & 6 & 4 \\
			6 & 4 & 5
		\end{matrix} \right)
		$$
		then $N^{(1)}$ is matrix of the bimodule $L$. Let $\textsf{C}^{(1)}$ be a bimodule category associated to $N^{(1)}$ with the above table of multiplication.
		
		Consider the matrix
		$$
		N^{(2)} = \left( \begin{matrix}
			\mathbb{Z}_3 & \emptyset \\
			R & \mathbb{Z}_3
		\end{matrix} \right)
		$$
		where
		$$
		\begin{matrix}
			\mathbb{Z}_3 = \{1,2,3\} & R = \{4,5,6\} & \mathbb{Z}_3 = \{7,8,9\}	
		\end{matrix} 
		$$
		such that
		$$
		\mathbb{Z}_3 \cdot R = \left( \begin{matrix}
			4 & 5 & 6 \\
			5 & 6 & 4 \\
			6 & 4 & 5
		\end{matrix} \right)
		$$
		and
		$$
		R \cdot \mathbb{Z}_3 = \left( \begin{matrix}
			4 & 6 & 5 \\
			5 & 4 & 6 \\
			6 & 5 & 4
		\end{matrix} \right).
		$$
		Similarly, $N^{(2)}$ is a matrix of the bimodule $R$. Notice that we changed the multiplication table of $R\cdot \mathbb{Z}_3$ by the involution of the group $\mathbb{Z}_3$. Let $\textsf{C}^{(2)}$ be a bimodule category associated to $N^{(2)}$ with the above table of multiplication.
		
		But the matrix
		$$
		N = \left( \begin{matrix}
			\mathbb{Z}_3 & L \\
			R & \mathbb{Z}_3
		\end{matrix} \right)
		$$
		with the above bimodules $\textsf{C}^{(1)}$ and $\textsf{C}^{(2)}$ doesn't admit a grouplike category. This was shown by calculating using \textsc{Mace4}. It also follows from Conclusion part (6) since the isomorphisms given by $L$ and $R$ are not inverses.
	\end{remark}

	We can conclude now the general structure of grouplike categories with only 2 objects $X_1, X_2$.
	\begin{center}
		$\begin{array}{c|ccc|ccc|ccc|ccc}
			&  & X_1 \rightarrow X_1 &  &  & X_1 \rightarrow X_2 &  &  & X_2 \rightarrow X_1 &  &  & X_2 \rightarrow X_2 &   \\ \hline
			X_1 &  &  &  &  &  & &  &  &  &  &  &  \\
			\downarrow &  & G^{*k_1} &  & & L &  &  &  &  &  &  &  \\ 
			X_1 &  &  &  &  &  & &  &  &  &  &  &  \\ \hline
			X_1 &  &  &  &  &  &  &  & &  &  &  &  \\
			\downarrow &  &  &  &  &  &  &  & G^{*i^{max}} &  &  & L & \\
			X_2 &  &  &  &  &  &  &  &  & &  &  & \\ \hline
			X_2 & &  &  &  & &  &  &  &  &  &  &  \\
			\downarrow &  & R &  &  & G^{*i^{max}} & &  &  &  &  &  &  \\
			X_1 & &  &  & & &  &  &  &  &  &  &  \\ \hline
			X_2 & &  &  &  &  &  &  &  &  &  & &  \\
			\downarrow &  &  &  &  &  &  & & R &  & & G^{*k_2} &  \\
			X_2 &  &  &  &  &  &  & &  & &  &  & 
		\end{array}$
	\end{center}
		
	\section{Applications}
	List all the 7 monoids of size 3:
	\begin{equation}\label{7monoids}
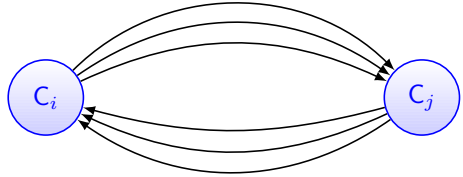
		\begin{array}{cccccccccc}
			\hline
			comp & & \textsf{C}_{1} & \textsf{C}_{2} & \textsf{C}_{3} & \textsf{C}_{4} & \textsf{C}_{5} & \textsf{C}_{6} &\textsf{C}_{7} \\
			\hline
			2^{2}    & = & 1 & 2 & 2 & 2 & 2 & 2 & 3 \\
			3^{2}    & = & 3 & 2 & 3 & 3 & 2 & 3 & 2 \\
			2\circ 3 & = & 3 & 2 & 2 & 3 & 3 & 2 & 1 \\
			3\circ 2 & = & 3 & 2 & 3 & 2 & 3 & 2 & 1 \\
			\hline
		\end{array}
	\end{equation}
	
	\begin{center}
		\begin{tikzpicture}[-latex ,auto ,node distance =4 cm and 5cm ,on grid ,
			semithick ,
			state/.style ={ circle ,top color =white , bottom color =blue!20 ,
				draw,blue , text=blue , minimum width =1 cm}]
			\node[state] (C)
			{$\textsf{C}_i$};
			\node[state] (A) [right=of C] {$\textsf{C}_j$};
			\path (C) edge [bend left =25] node[below =0.10 cm] {} (A);
			\path (C) edge [bend left =35] node[below =0.001 cm] {} (A);
			\path (C) edge [bend left =45] node[below =0.002 cm] {} (A);
			\path (A) edge [bend right = -15] node[below =0.15 cm] {} (C);
			\path (A) edge [bend right = -25] node[below =0.15 cm] {} (C);
			\path (A) edge [bend right = -35] node[below =0.15 cm] {} (C);
		\end{tikzpicture}
	\end{center}
	We combine monoids of 3 elements together in one category, two monoids $\textsf{C}_i, \textsf{C}_j$ are called \textit{connected} if there exists a category where the monoids $\textsf{C}_i$ and $\textsf{C}_j$ are the endomorphism monoids of the two objects. Viewing them as objects, each object is one of the monoid of endomorphisms listed before, and this graph of a category is associated to the matrix:
	\begin{center}
		$\left( \begin{matrix}
			3 & 3 \\
			3 & 3
		\end{matrix} \right)$.
	\end{center}
	
	From Table \ref{7monoids}, we can remark that $\textsf{C}_5$ and $\textsf{C}_6$ are grouplike monoids. Where $\textsf{C}_5$ is of the form $G^{*1}$ such that $G= \mathbb{Z}_2$, and $\textsf{C}_6$ is of the form $G^{*2}$ such that $G = \{2\}$ the trivial group.

	\begin{example} (\textbf{The number of categories between $\textsf{C}_6$ and itself})\\
		Let $\textsf{B}$ be a category with two objects $X$ and $Y$ associated to the matrix
		$$
		M_1 = \left( \begin{array}{cc}
			3 & 3 \\
			0 & 3
		\end{array} \right)
		$$
		such that $\textsf{C}(X,X) = \textsf{C}_6$, $\textsf{C}(Y,Y) = \textsf{C}_6$, $\textsf{C}(X,Y) = L$ and $\textsf{C}(Y,X) = \emptyset$. Fix
		$$
		\textsf{C}_6 = \{1,2,3\} ~~~~~ L = \{4,5,6\} ~~~~~ \textsf{C}_6 = \{7,8,9\}
		$$ 
		(these are two distinct copies of $\textsf{C}_6$).
		\begin{enumerate}
			\item $\textbf{i = 0:}$ the number of categories associated to $M_1$ is 64 (calculated by \textsc{Mace4} \cite{mace4}). But Since $\textsf{C}_6$ is a grouplike category containing the trivial group $\{2\}$, then the orbit of $L$ has size 1 and it's unique. This means that $2 * x$ is unique for all $x \in L$. Suppose it's equal to 4. This reduces the number of possibilities to 15.
			
			Similarly, let $\textsf{D}$ be a category with two objects $X$ and $Y$ associated to the matrix
			$$
			M_2 = \left( \begin{array}{cc}
				3 & 0 \\
				3 & 3
			\end{array} \right)
			$$
			such that $\textsf{C}(X,X) = \textsf{C}_6$, $\textsf{C}(Y,Y) = \textsf{C}_6$, $\textsf{C}(X,Y) = \emptyset$ and $\textsf{C}(Y,X) = R$.
			
			For the same reason, the number of possible categories associated to $M_2$ with the orbit condition is 15.
			
			Then a category $\textsf{C}$ associated to the matrix
			$$
			M = \left( \begin{array}{cc}
				3 & 3 \\
				3 & 3
			\end{array} \right)
			$$
			such that $\textsf{C}(X,X) = \textsf{C}_6$, $\textsf{C}(Y,Y) = \textsf{C}_6$, $\textsf{C}(X,Y) = L$ and $\textsf{C}(Y,X) = R$, has 15 left bimodules $A_1,\hdots,A_{15}$ and 15 right bimodules $B_1,\hdots,B_{15}$.
			
			We have two cases:
			\begin{enumerate}
				\item The bimodules are the same
				$$
				\left(\begin{array}{cc}
					G^{*2} & A_i \\
					B_i & G^{*2}
				\end{array}
				\right).
				$$
				By Theorem \ref{grouplikecounting}, there are exactly 15 possibilities for this matrix (up to isomorphism).
				\item The bimodules are different
				$$
				\left(\begin{array}{cc}
					G^{*2} & A_i \\
					B_j & G^{*2}
				\end{array}
				\right)
				$$
				such that $1 \leq i < j \leq 15$.
				
				\begin{center}
					
					\begin{tabular}{|l|*{10}{c}}\hline
						\backslashbox{}{$\frac{15 \times 14}{2}$} \\
						\hline
					\end{tabular}
				\end{center}
				There are $\frac{15 \times 14}{2} = 105$ possibilities.
			\end{enumerate}
			Hence for $i = 0$, we have 120 categories.
			
			\item $\textbf{i = 1:}$ there are 2 left bimodules $A_1,A_2$ and 2 right bimodules $B_1,B_2$. Then the possibilities are:
			\begin{enumerate}
				\item The bimodules are the same
				$$
				\left(\begin{array}{cc}
					G^{*2} & A_i \\
					B_i & G^{*2}
				\end{array}
				\right)
				$$
				and we have 2 categories that could be associated to this matrix.
				\item The bimodules are different
				$$
				\left(\begin{array}{cc}
					G^{*2} & A_i \\
					B_j & G^{*2}
				\end{array}
				\right)
				$$
				and we have 1 category (up to isomorphism).
			\end{enumerate}
			Hence in total there are $120 + 3 = 123$ categories between $\textsf{C}_6$ and itself.
		\end{enumerate}
	\end{example}
	
	\begin{example}\textbf{(The number of categories between $\textsf{C}_5$ and itself)}\\
		The number of categories between $\textsf{C}_5$ and itself can also be found bimodules. Let $\textsf{B}$ be a category with two objects $X$ and $Y$ associated to the matrix 
		$$
		M_1 = \left( \begin{array}{cc}
			3 & 3 \\
			0 & 3
		\end{array} \right)
		$$
		such that $\textsf{C}(X,X) = \textsf{C}_5$, $\textsf{C}(Y,Y) = \textsf{C}_5$, $\textsf{C}(X,Y) = L$ and $\textsf{C}(Y,X) = \emptyset$.
		
		And similarly let $\textsf{D}$ be a category with two objects $X$ and $Y$ associated to the matrix 
		$$
		M_2 = \left( \begin{array}{cc}
			3 & 0 \\
			3 & 3
		\end{array} \right)
		$$
		such that $\textsf{C}(X,X) = \textsf{C}_5$, $\textsf{C}(Y,Y) = \textsf{C}_5$, $\textsf{C}(X,Y) = \emptyset$ and $\textsf{C}(Y,X) = R$.
		
		Then there exists a category $\textsf{C}$ with two objects associated to the matrix
		$$
		M = \left( \begin{array}{cc}
			3 & 3 \\
			3 & 3
		\end{array} \right)
		$$
		such that $\textsf{C}(X,X) = \textsf{C}_5$, $\textsf{C}(Y,Y) = \textsf{C}_5$, $\textsf{C}(X,Y) = L$ and $\textsf{C}(Y,X) = R$. Where $\textsf{B}$ is a left bimodule and $\textsf{D}$ is a right bimodule.
		
		Since there is 1 left and 1 right bimodule (also calculated by \textsc{Mace4}) and since $i^{max} =0$, then the number of possible $y_0$ is equal to the order of the center of the group $\mathbb{Z}_2$ which is 2 (Theorem \ref{grouplikecounting} and Remark \ref{center}). Hence there are 2 categories between $\textsf{C}_5$ and itself.
	\end{example}
	
	\begin{definition}
		A category $\textsf{C}$ is called \textit{reduced} if there does not exist any two distinct isomorphic objects in $\textsf{C}$.
	\end{definition}
	
	\begin{lemma} \label{lemothermehtod}
		Let $\textsf{C}$ be a reduced category associated to the matrix
		$$
		\left(\begin{array}{cc}
			A & L \\
			R & B
		\end{array}\right).
		$$
		If there exists $x \in L$ and $y \in R$ such that $y \cdot x = 1_{B}$, then $|B| < |A|$ and $B$ is a sub-monoid of $A$ disjoint from $\{1_{A}\}$.
		\begin{proof}
			Let
			\begin{center}
				$	\begin{array}{ccccc}
					\phi & : & A & \rightarrow & B \\
					&	 & b & \mapsto & x\cdot f \cdot y
				\end{array} $
			\end{center}
			We have
			$$
			\phi(ff') = x \cdot ff' \cdot y = x \cdot f \cdot 1_B \cdot f' \cdot y = (x \cdot f \cdot y) \cdot (x \cdot f' \cdot y) = \phi(f)\phi(f').
			$$
			$\phi$ is injective, indeed, suppose that $\phi(f) = 1_{A}$, then
			\begin{eqnarray*}
				x\cdot f\cdot y &=& 1_{A} \\
				\Rightarrow	x \cdot f \cdot y \cdot x &=& x \\
				\Rightarrow	x \cdot f &=& x \\
				\Rightarrow	y \cdot x \cdot f &=& y \cdot x\\
				\Rightarrow	f &=& 1_{B}.
			\end{eqnarray*}
		\end{proof}
	\end{lemma}
	
	\begin{remark}
	We note that if $|A| = |B|$ then $x \cdot y = 1_{A}$ and $\textsf{C}$ is not reduced.
	\end{remark}
	
	\begin{prop} \label{grouplikeconnected}
		Let $M = \left( \begin{matrix}
			3 & a \\
			b & 3 
		\end{matrix} \right) $ and $ \textsf{C}$ be a category associated to $M$ whose objects are $X$ and $Y$. Then
		$$
		\textsf{C}(Y,Y) = \textsf{C}_5 \iff \textsf{C}(X,X)= \textsf{C}_5.
		$$
		\begin{proof}
		Suppose that $\textsf{C}(X,X) = A, \textsf{C}(X,Y) = L, \textsf{C}(Y,X) = R$ and $\textsf{C}(Y,Y) = \textsf{C}_5 = \mathbb{Z}_2 \cup \{1\}$. Consider the algebraic matrix of $\textsf{C}$
		$$
		M^{alg} = \left(\begin{array}{cc}
		A & L \\
		R & \textsf{C}_5
		\end{array}\right)
		$$
		such that $|A| = 3, |L| = a, |R| = b$ and $|\textsf{C}_5| = 3$.
		
		Let $\textsf{C}'$ be a subcategory of $\textsf{C}$ defined in the following way
		\begin{enumerate}
			\item \textbf{(Objects of $\textsf{C}'$):} $X$ and $Y$. 
			\item \textbf{(Morphisms of $\textsf{C}'$):} $A, L\cdot \mathbb{Z}_2, \mathbb{Z}_2 \cdot R$ and $\mathbb{Z}_2$. (It means that we take the orbits of $L$ and $R$ by the action of $\mathbb{Z}_2$).
			\item \textbf{(Composition):} by Proposition \ref{compoforbits} we have that the composition of orbits goes into the group $\mathbb{Z}_2$.
			\item \textbf{(Identities):} $1_A$ and $1_{\mathbb{Z}_2}$.
		\end{enumerate}
		We obtain that $\textsf{C}'$ is associated to the following algebraic matrix 
		$$
		M'^{alg} = \left( \begin{matrix}
			A & L\cdot \mathbb{Z}_2 \\
			\mathbb{Z}_2 \cdot R & \mathbb{Z}_2 
		\end{matrix} \right).
		$$
		
		Since there exist $x \in L \cdot \mathbb{Z}_2$ and $y \in \mathbb{Z}_2 \cdot R$ such that $x \cdot y = 1_{\mathbb{Z}_2}$, then from Lemma \ref{lemothermehtod} we have that $\mathbb{Z}_2$ is a sub-monoid of $A$ disjoint from $\{1_{A}\}$. Therefore, $A = \mathbb{Z}_2 \cup \{1\} = \textsf{C}_5$.	
		\end{proof}
	\end{prop}
	
	\begin{theorem}
		Let $M = \left( \begin{matrix}
			3 & a_{12}& \hdots & a_{1n}\\
			a_{21} & 3 & \hdots & a_{2n} \\
			\vdots &\vdots & \ddots & \vdots \\
			a_{n1} & \hdots & \hdots & 3
		\end{matrix} \right)$ be a strictly positive matrix and let $\textsf{C}$ be a reduced category associated to $M$ with $Ob(\textsf{C}) = \{1,\hdots ,n\}$.\\
		If there exists $i\in Ob(\textsf{C})$ such that $\textsf{C}(i,i) = \textsf{C}_5$ then $\textsf{C}(j,j) = \textsf{C}_5~~ \forall j \in Ob(\textsf{C})$.
		\begin{proof}
			If $\textsf{C}$ is a category associated to $M$ then every regular sub-matrix of $M$ of size 2 is associated to a full sub-category of $\textsf{C}$ \cite{samer1} with two objects, then we apply Proposition \ref{grouplikeconnected} to get the result.
		\end{proof}
	\end{theorem}
	
	The classification of finite categories depends mostly on the algebraic structure of the endomorphism monoids. Monoids that are groups or contain a group impose some restrictions on the cardinality of the sets of morphisms. The properties obtained by studying such monoids makes the classification and the counting problem easier. Counting finite structures is not easy in general, but optimizing the number of categories that could be obtained is a good start. We are now studying other types of monoids hoping to get a similar classification theorems as in this paper.
	
	\section*{Acknowledgment}
	I would like to thank Wesley Fussner for he is the one that suggested working on grouplike monoids and categories, and Angel Toledo for the reference to Leinster's paper. A very big thanks to my advisor Carlos Simpson for his constant help during my thesis.

	\bibliographystyle{siam}
	\bibliography{biblio}

\end{document}